\documentclass[a4paper,12pt,reqno]{amsart}
\usepackage{amssymb}
\usepackage{amsmath}
\usepackage{amsthm}
\usepackage{amsfonts}
\usepackage{mathrsfs}

\newtheorem{theorem}{Theorem}[section]

\newtheorem{prop}[theorem]{Proposition}

\theoremstyle{definition}
\newtheorem{definition}[theorem]{Definition}

\theoremstyle{remark}
\newtheorem{remark}[theorem]{Remark}

\newcommand{\ra}{\rightarrow}

\numberwithin{equation}{section}

\allowdisplaybreaks[1]

\begin{document}

\title[R-duals for Gabor system]{On Riesz duals for the Gabor system on LCA groups}

\author{S. Arati*}

\author{P. Devaraj}

\address{School of Mathematics, Indian Institute of Science Education and Research Thiruvananthapuram, Thiruvananthapuram 695551, India}

\email{aratishashi@iisertvm.ac.in, devarajp@iisertvm.ac.in}
\thanks{*Corresponding author}



\begin{abstract}
In this paper, we analyse the circumstances in which the adjoint Gabor system is an R-dual of a given Gabor frame in the context of separable uniform time-frequency lattices in locally compact abelian groups. In this regard, we also prove a necessary condition for a Gabor Bessel sequence in this setting to be complete.
\end{abstract}

\maketitle
\section{Introduction}
The study of Gabor theory, which is a widely researched topic, gained popularity not only in the classical case of the real line but also in the abstract context, namely the locally compact abelian group setting. Let $G$ be a second countable locally compact abelian (LCA) group and $\widehat{G}$ be its dual group, which consists of all the continuous characters of $G$. For $u\in G$ and $\xi\in\widehat{G}$, the translation $\mathcal{T}_{u}$ and modulation $E_{\xi}$ operators on $L^2(G)$ are given by
\begin{align*}
\mathcal{T}_{u}f(x)=f(x-u)\text{ and }E_{\xi}f(x)=\xi(x)f(x),\quad \text{for }f\in L^2(G).
\end{align*}
A Gabor system in groups is generated using these translations and modulations. Back in 1998, Gr{\"o}chenig in \cite{Gro1998aspects} showed that not all results in Gabor theory on $\mathbb{R}$ can be easily extended to the group setting analogously. For instance, even the weaker form of the Balian-Low theorem, in terms of the short time Fourier transform, does not hold true in general for all LCA groups. It has been proved to be false for discrete and compact groups. However, in \cite{Kani}, a variant of it has been shown to be true for compactly generated LCA groups under certain conditions.
Recently, in 2020, Enstad in \cite{Enstad2020balian} has established a Balian-Low theorem, in the case of a non-compactly generated group, namely $\mathbb{R}\times\mathbb{Q}_p$, where $\mathbb{Q}_p$ represents the p-adic numbers.
\par
Two amongst the several fundamental principles in Gabor analysis are the density and the duality theorems. In the context of locally compact abelian groups, a critically sampled Gabor system is characterized to be a frame using the Zak transform in groups in \cite{Gro1998aspects}. Also, sufficient conditions for when a Gabor system, generated by separable uniform lattices, forms a frame as well as necessary conditions for the same, in terms of the density, are given. The duality and biorthogonality relations in the case of non-separable uniform lattices for elementary LCA groups have been studied by Feichtinger and Kozek in \cite{feich1998quant}. In \cite{jakobsen2016co}, Jakobsen and Lemvig have considered co-compact Gabor systems, which are Gabor systems generated by separable, closed, co-compact  subgroups (not necessarily discrete) and proved several results pertaining to Gabor analysis, including the Walnut's and Janssen's representations of the frame operator, the Wexler-Raz biorthogonality relations and the duality principle. The same authors in \cite{jakobsen2016density} have proved density theorems in the case of Gabor systems generated using time-frequency shifts along non-separable, closed subgroups of $G\times\widehat{G}$, wherein the necessary conditions are provided in terms of a generalization of the density measure, which is the measure of a fundamental domain in the case of uniform lattices. Besides, they have also dealt with duality principles and the Wexler-Raz biorthogonality relations in this context. Recently in 2021, Jakobsen and Luef have explored the duality between multi-window super Gabor systems along with Heisenberg modules in \cite{jakobsen2021duality}. The pioneer work on the duality principle for Gabor systems on the real line is due to Daubechies, Landau, Landau\cite{daub1994gabor}, Janssen\cite{janssen1994duality} and Ron, Shen\cite{ron1997weyl}.
A detailed account of the density and duality theorems has been provided for various scenarios on the reals, such as rectangular lattices in one dimension and higher dimensions as well as irregular Gabor systems, in \cite{heil2007history} and the references therein.
\par
Interestingly, Casazza, Kutyniok and Lammers, in \cite{casazza2004duality}, have considered the problem of extending the duality principles and the Wexler-Raz biorthogonality relations in Gabor analysis on the reals to the arbitrary frame setting by employing the concept of Riesz duals(or R-duals) introduced by them. It turns out that the duality relation between a given frame and a corresponding R-dual is the same as that between a Gabor frame for $L^2(\mathbb{R})$ and its adjoint Gabor system, which naturally raises the question if the adjoint Gabor system is indeed an R-dual of the given Gabor frame. This has been proved to be true in the case when the given Gabor frame is a Riesz basis or a tight frame for $L^2(\mathbb{R})$. For the other cases, it remains as an open problem. Later, in 2016, Stoeva and Christensen\cite{stoeva2016various} altered the concept of R-duals and defined the modified R-duals of type III, which in fact generalizes the duality principle in Gabor theory while at the same time retaining the optimality of the bounds. For characterizations of the different types of R-duals, one may refer to \cite{stoeva2015r, chuang2015equiv}. Moreover, the study on R-duals has engendered further related research along various avenues. See \cite{xiao2009duality, chrisXiao2013charact, enayati2017duality, takhteh2017r, li2019charact, dong2020duality, li2020class, li2021duality, li2021g-duality} in this regard.
\par
The connection between a Gabor frame and its adjoint Gabor system in the locally compact abelian group setting is also similar to that between a frame and its R-dual, which has been established in  \cite{jakobsen2016co} (see Theorem \ref{T:duality in LCA}.). Hence, the question of whether the R-duals will generalize the duality principle in Gabor theory even in this abstract context is an interesting one, which will be addressed in this paper. For a second countable infinite locally compact abelian (LCA) group $G$ and its dual group $\widehat{G}$, let $\Lambda$ and $\Gamma$ respectively denote countable uniform lattices. In other words, they are countable discrete subgroups, which are co-compact. These uniform lattices generalize the full-rank lattices in higher dimension of reals. Further, let $\Lambda^{\perp}\subset \widehat{G}$ and $\Gamma^{\perp}\subset G$ be the annihilators of $\Lambda$ and $\Gamma$ respectively. For a given $g\in L^2(G)$, we are interested in the existence of orthonormal bases $\{e_{\alpha\beta}\}_{(\alpha,\beta)\in\Gamma^{\perp}\times\Lambda^{\perp}}$ and $\{h_{\lambda\gamma}\}_{(\lambda,\gamma)\in\Lambda\times\Gamma}$ for the separable Hilbert space $L^2(G)$, which is meaningful to expect as both $\Lambda\times\Gamma$ as well as $\Gamma^{\perp}\times\Lambda^{\perp}$ are infinite in the case of an infinite group $G$, such that 
\[E_{\beta}\mathcal{T}_\alpha g=\sum_{(\lambda,\gamma)\in\Lambda\times\Gamma}\langle E_{\gamma}\mathcal{T}_{\lambda}g, e_{\alpha\beta}\rangle_{L^2(G)}h_{\lambda\gamma},\quad(\alpha,\beta)\in\Gamma^{\perp}\times\Lambda^{\perp}.\]
In other words, this would mean to ask if the adjoint Gabor system $\{E_{\beta}\mathcal{T}_\alpha g\}_{(\alpha,\beta)\in\Gamma^{\perp}\times\Lambda^{\perp}}$ is an R-dual of the given Gabor system $\{E_{\gamma}\mathcal{T}_{\lambda}g\}_{(\lambda,\gamma)\in\Lambda\times\Gamma}$. We shall answer this affirmatively in certain circumstances in Section \ref{S:Adjoint Gabor system as R-dual}. More precisely, we first show that under certain conditions on the size of the time-frequency lattice, the adjoint Gabor system $\{E_{\beta}\mathcal{T}_\alpha f\}_{(\alpha,\beta)\in\Gamma^{\perp}\times\Lambda^{\perp}}$ is an R-dual, when the given Gabor frame $\{E_{\gamma}\mathcal{T}_{\lambda}f\}_{(\lambda,\gamma)\in\Lambda\times\Gamma}$ is tight. Interestingly, it also turns out that the adjoint Gabor system is an R-dual of any critically sampled Gabor system. In fact, this result is true in the context of finite abelian groups as well. In connection with these R-duality results, we also prove a necessary condition when the given Gabor system is a complete Bessel sequence, which is discussed in Section \ref{S:completeness}. That the completeness and frame property of a Gabor system are significantly different is explicated in \cite{grochenig2016complete}.
\par
We also wish to remark that the question on the R-duality of the adjoint Gabor system in all the other cases remains as an open problem.

\section{Necessary background}

We shall now recall some preliminaries needed from frame theory. Let $0\neq\mathcal{H}$ be a separable Hilbert space. A sequence $\{f_{k}:k\in\mathbb{N}\}$ in $\mathcal{H}$ is a frame for $\mathcal{H}$ if there exist constants $A,B>0$ such that \[A\|f\|^{2}\leq\sum_{k\in\mathbb{N}}|\langle f,f_k\rangle|^{2}\leq B\|f\|^{2},\quad \forall\, f\in\mathcal{H}.\]The constants $A$ and $B$ are called the frame bounds. If the right hand side inequality holds, then $\{f_{k}:k\in\mathbb{N}\}$ is said to be a Bessel sequence with bound $B$. A frame is called a tight frame if $A=B$ and a Parseval frame if $A=B=1$ in the above inequalities. A sequence $\{f_k:k\in\mathbb{N}\}$ in $\mathcal{H}$ is said to be a frame sequence if it is a frame for $\overline{\text{span}}\{f_k:k\in\mathbb{N}\}$. A Riesz basis for $\mathcal{H}$ is a family of the form $\{Ue_{k}:k\in\mathbb{N}\}$, where $\{e_{k}:k\in\mathbb{N}\}$ is an orthonormal basis for $\mathcal{H}$ and $U:\mathcal{H}\rightarrow\mathcal{H}$ is a bounded invertible operator. Alternatively, a sequence $\{f_{k}:k\in\mathbb{N}\}$ is a Riesz basis for $\mathcal{H}$ if $\{f_{k}:k\in\mathbb{N}\}$ is complete in $\mathcal{H}$, and there exist constants $A,B>0$ such that for every finite scalar sequence $\{c_k\}$, one has
\[A\sum_{k}|c_k|^2\leq \left\|\sum_{k}c_k f_k\right\|^2\leq B\sum_{k}|c_k|^2.\]
A sequence $\{f_k:k\in\mathbb{N}\}$ in $\mathcal{H}$ is called a Riesz sequence if it is a Riesz basis  for $\overline{\text{span}}\{f_k:k\in\mathbb{N}\}$.
\par
For a Bessel sequence $F=\{f_{k}:k\in\mathbb{N}\}$ in $\mathcal{H}$, the operator $T_{F}:\mathcal{H}\ra l^2(\mathbb{N})$ given by $T_{F}(g)=\{\langle g,f_{k} \rangle\}_{k\in\mathbb{N}}$, $g\in\mathcal{H}$ is the analysis operator and its adjoint operator $T^{\ast}_{F}:l^2(\mathbb{N})\ra \mathcal{H}$, given by $T^{\ast}_{F}(\{c_{k}\})=\sum\limits_{k\in \mathbb{N}}c_{k}f_{k}$ is the synthesis operator. For a sequence $F=\{f_k:k\in\mathbb{N}\}$ in $\mathcal{H}$, the deficit $\mathcal{D}(F)$ is defined to be the least number of elements which must be adjoined in order to make it complete, while the excess $\mathcal{E}(F)$ is defined to be greatest number of elements that can be deleted without losing completeness. More formally,
\begin{align*}
\mathcal{D}(F)&=\inf\big\{|X|:X\subset\mathcal{H},\,\overline{span}(F\cup X)=\mathcal{H}\big\}\\
\text{and }\qquad
\mathcal{E}(F)&=\sup\big\{|X|:X\subset F,\,\overline{span}(F\setminus X)=\overline{span}(F)\big\}.
\end{align*}
It is proved in \cite{Radu} that the deficit and excess of a Bessel sequence $F=\{f_k:k\in\mathbb{N}\}$ are related with the analysis and synthesis operators as follows.
\begin{align}
\text{(i) }&\mathcal{D}(F)=dim(ker(T_{F}))\text{ and }\mathcal{E}(F)\geq dim(ker(T^{\ast}_{F})).\nonumber\\
\text{(ii) }&\text{If }F\text{ is a frame, then }\mathcal{E}(F)=dim(ker(T^{\ast}_{F})).\label{E:excess and synth ker}
\end{align}
More is said, using the concepts of intertwining operators and point spectrum. A pair of bounded operators $M:\mathcal{H}\ra \mathcal{H}$ and $N:l^2(\mathbb{N})\ra l^2(\mathbb{N})$ satisfying $NT_{F}=T_{F}M$ is called an intertwining pair of operators for $F=\{f_k\}_{k\in \mathbb{N}}$. We recall that an operator $T$ on $\mathcal{H}$ has no point spectrum if there are no complex values $a$ such that ker$(T-aI)\neq\{0\}$.
\begin{theorem}\label{T:intertwining}
Suppose $(M,N)$ is an intertwining pair of operators for a Bessel sequence $F=\{f_k\}_{k\in \mathbb{N}}$ in $\mathcal{H}$. Then,
\begin{enumerate}
\item[(i)] if $M$ has no point spectrum, then $\mathcal{D}(F)$ is either zero or infinity.
\item[(ii)] if $N^{\ast}$ has no point spectrum and $F$ is a frame, then $\mathcal{E}(F)$ is either zero or infinity.
\end{enumerate}
\end{theorem}

\par

The following are a few results related to the frame properties of the Gabor system $\{E_{\gamma}\mathcal{T}_{\lambda}g\}_{(\lambda,\gamma)\in\Lambda\times\Gamma}$ and the adjoint Gabor system $\{E_{\beta}\mathcal{T}_\alpha g\}_{(\alpha,\beta)\in\Gamma^{\perp}\times\Lambda^{\perp}}$ in $L^2(G)$.
\begin{theorem}[\cite{jakobsen2016co}]
The Gabor system $\{E_{\gamma}\mathcal{T}_{\lambda}g\}_{(\lambda,\gamma)\in\Lambda\times\Gamma}$ is a Bessel sequence with bound $B>0$ if and only if $\{E_{\beta}\mathcal{T}_\alpha g\}_{(\alpha,\beta)\in\Gamma^{\perp}\times\Lambda^{\perp}}$ is a Bessel sequence with bound $B$.
\end{theorem}

\begin{theorem}[\cite{jakobsen2016co}]\label{T:duality in LCA}
The Gabor system $\{E_{\gamma}\mathcal{T}_{\lambda}g\}_{(\lambda,\gamma)\in\Lambda\times\Gamma}$ is a frame for $L^2(G)$ with bounds $A$ and $B$ if and only if $\{E_{\beta}\mathcal{T}_\alpha g\}_{(\alpha,\beta)\in\Gamma^{\perp}\times\Lambda^{\perp}}$ is a Riesz sequence with bounds $A$ and $B$.
\end{theorem}

\begin{theorem}[\cite{jakobsen2016co}]\label{T:dual orthog}
The Gabor system $\{E_{\gamma}\mathcal{T}_{\lambda}g\}_{(\lambda,\gamma)\in\Lambda\times\Gamma}$ is a tight frame for $L^2(G)$ if and only if $\{E_{\beta}\mathcal{T}_\alpha g\}_{(\alpha,\beta)\in\Gamma^{\perp}\times\Lambda^{\perp}}$ is an orthogonal system. In this case, the frame bound $A=\|g\|^2_{L^2(G)}$.
\end{theorem}

The concept of R-duals and the related duality principle in the abstract frame theory, which turns out to be similar to the result given in Theorem \ref{T:duality in LCA}, are as follows.

\begin{definition}[\cite{casazza2004duality}]
Let $\{e_j\}_{j\in\mathbb{N}}$ and $\{h_i\}_{i\in\mathbb{N}}$ be orthonormal bases for a separable Hilbert space $\mathcal{H}$. Let $\{f_i\}_{i\in\mathbb{N}}\subset\mathcal{H}$ be such that $\sum\limits_i|\langle f_i,e_j\rangle|^2<\infty$ for all $j\in\mathbb{N}$, and let $w_j^f=\sum\limits_{i=1}^\infty\langle f_i,e_j\rangle h_i,\; j\in\mathbb{N}$.
Then, $\{w_j^f\}$ is called the R-dual sequence for $\{f_i\}$ with respect to $\{e_j\}$ and $\{h_i\}$.
\end{definition}

\begin{theorem}[Ron-Shen duality principle in abstract frame theory-\cite{casazza2004duality}]
Suppose there exist positive constants $A, B$ such that $A\leq\sum\limits_i|\langle f_i,e_j\rangle|^2\leq B, \, j\in\mathbb{N}$. Then, $\{f_i\}$ is a frame for $\mathcal{H}$ with frame bounds $A$ and $B$ if and only if $\{w_j^f\}$ is a Riesz sequence with Riesz bounds $A$ and $B$.
\end{theorem}

We answer the question of the above duality principle extending Theorem \ref{T:duality in LCA}, based on the size of the lattice, $d(\Lambda\times\Gamma)$, given by 
\begin{align*}
d(\Lambda\times\Gamma)&:=\mu_{(G\times\widehat{G})/(\Lambda\times\Gamma)}
\left((G\times\widehat{G})/(\Lambda\times\Gamma)\right)
=\mu_{G/\Lambda}\left(G/\Lambda\right)\mu_{\widehat{G}/\Gamma}(\widehat{G}/\Gamma)\\
&=\mu_G(\mathcal{S}_{\Lambda})\mu_{\widehat{G}}(\mathcal{S}_{\Gamma}),
\end{align*}
where $\mu_G$ denotes the Haar measure of the LCA group $G$ and $\mathcal{S}_{\Lambda}$ denotes a Borel measurable section of $G/\Lambda$, the existence of which is well-known. The reciprocal of $d(\Lambda\times\Gamma)$ gives the density measure of the lattice. It is proved in \cite{Gro1998aspects} that $d(\Gamma^{\perp}\times\Lambda^{\perp})=\dfrac{1}{d(\Lambda\times\Gamma)}$. The following is an important and useful identity in abstract group theory.

\begin{definition}[Weil's formula]
Let $G$ be a LCA group and $H$ be its closed subgroup. Then, 
\[\int_G f(x)d\mu_G(x)=\int_{G/H}\int_{H}f(x+y)d\mu_H(y)d\mu_{G/H}(x+H),\quad f\in L^1(G),\]
with the constant factor in the Haar measure $\mu_{G/H}$ on the quotient group $G/H$ being appropriately chosen.
\end{definition}
For a detailed study on abstract group theory, one may refer to \cite{rudin2017fourier, Foll2}.

\section{Completeness of the Gabor system}\label{S:completeness}
The following theorem provides a necessary condition for the Gabor system $\{E_{\gamma}\mathcal{T}_{\lambda}f\}_{(\lambda,\gamma)\in\Lambda\times\Gamma}$ to be a complete Bessel sequence, in terms of the size of the time-frequency lattice $\Lambda\times\Gamma$.
\begin{theorem}\label{T:GF completeness}
Suppose there exist a finite number, say $k$, of elements $\alpha_i\in \Gamma^{\perp}$ such that $\bigcup\limits_{i=1}^k (\alpha_i+\mathcal{S}_{\Gamma^{\perp}})\cap\mathcal{S}_{\Lambda}=\mathcal{S}_{\Lambda}$. For a function $f\in L^2(G)$, if the Gabor system $\{E_{\gamma}\mathcal{T}_{\lambda}f\}_{(\lambda,\gamma)\in\Lambda\times\Gamma}$ is a complete Bessel sequence, then $d(\Lambda\times\Gamma)\leq 1$.
\end{theorem}

\begin{proof}
Let  $\{E_{\gamma}\mathcal{T}_{\lambda}f\}_{(\lambda,\gamma)\in\Lambda\times\Gamma}$, $f\in L^2(G)$ be a complete Bessel sequence. Consider the operator $T^{\ast}_{f:\Lambda,\Gamma}T_{f:\Lambda,\Gamma}:L^2(G)\ra L^2(G)$, where $T_{f:\Lambda,\Gamma}$ and $T^{\ast}_{f:\Lambda,\Gamma}$ denote the analysis and synthesis operators respectively. More explicitly, the analysis operator $T_{f:\Lambda,\Gamma}:L^2(G)\ra l^2(\Lambda\times\Gamma)$ is given by $T_{f:\Lambda,\Gamma}(g)=\{\langle g,E_{\gamma}\mathcal{T}_{\lambda}f \rangle\}_{(\lambda,\gamma)\in\Lambda\times\Gamma}$ and its adjoint  $T^{\ast}_{f:\Lambda,\Gamma}:l^2(\Lambda\times\Gamma)\ra L^2(G)$, the synthesis operator, is given by $T^{\ast}_{f:\Lambda,\Gamma}(\{c_{\lambda\gamma}\})=\sum\limits_{(\lambda,\gamma)\in\Lambda\times\Gamma}c_{\lambda\gamma}E_{\gamma}\mathcal{T}_{\lambda}f$. As $\{E_{\gamma}\mathcal{T}_{\lambda}f\}_{(\lambda,\gamma)\in\Lambda\times\Gamma}$ is a Bessel sequence, the operator $T^{\ast}_{f:\Lambda,\Gamma}T_{f:\Lambda,\Gamma}$ is bounded, self-adjoint and positive. However, it need not be invertible. Hence, we consider the operator $R_{\theta}:=\theta I+T^{\ast}_{f:\Lambda,\Gamma}T_{f:\Lambda,\Gamma}$ for $\theta>0$. This operator is in fact bounded, self-adjoint and bounded below, and thus, invertible. Let $R_{\theta}^{-1}(f)$ be $h_{\theta}$.
\par
Consider, for $g\in L^2(G)$,
\begin{align*}
T_{f:\Lambda,\Gamma}R_{\theta}^{-1}(g)
&=\{\langle R_{\theta}^{-1}g,E_{\gamma}\mathcal{T}_{\lambda}f \rangle\}_{(\lambda,\gamma)\in\Lambda\times\Gamma}\\
&=\{\langle g, R_{\theta}^{-1}E_{\gamma}\mathcal{T}_{\lambda}f \rangle\}_{(\lambda,\gamma)\in\Lambda\times\Gamma}.
\end{align*}
The operators $E_{\gamma}$ and $\mathcal{T}_{\lambda}$ satisfy the following commutativity relation
\begin{equation}\label{E:commu relation}
E_{\gamma}\mathcal{T}_{\lambda}=\gamma(\lambda)\mathcal{T}_{\lambda}E_{\gamma}\text{ for }\lambda\in\Lambda,\,\gamma\in \Gamma.
\end{equation}
Now, the operator $T^{\ast}_{f:\Lambda,\Gamma}T_{f:\Lambda,\Gamma}$ commutes with $E_{\gamma}\mathcal{T}_{\lambda}$, $\lambda\in\Lambda,\, \gamma\in\Gamma$ for, taking $g\in L^2(G)$, we have
\begin{align*}
T^{\ast}_{f:\Lambda,\Gamma}T_{f:\Lambda,\Gamma}E_{\gamma}\mathcal{T}_{\lambda}(g)
&=\sum_{\substack{\lambda'\in\Lambda\\\gamma'\in\Gamma}}\langle E_{\gamma}\mathcal{T}_{\lambda}g, E_{\gamma'}\mathcal{T}_{\lambda'}f\rangle E_{\gamma'}\mathcal{T}_{\lambda'}f\\
&=\sum_{\substack{\lambda'\in\Lambda\\\gamma'\in\Gamma}}\langle g, \mathcal{T}_{-\lambda}E_{\gamma'-\gamma}\mathcal{T}_{\lambda'}f\rangle E_{\gamma'}\mathcal{T}_{\lambda'}f\\
&=\sum_{\substack{\lambda'\in\Lambda\\\gamma'\in\Gamma}}\langle g, (\gamma'-\gamma)(\lambda)E_{\gamma'-\gamma}\mathcal{T}_{\lambda'-\lambda}f\rangle E_{\gamma'}\mathcal{T}_{\lambda'}f\\
&=\sum_{\substack{\lambda'\in\Lambda\\\gamma'\in\Gamma}}\langle g, E_{\gamma'-\gamma}\mathcal{T}_{\lambda'-\lambda}f\rangle \overline{(\gamma'-\gamma)(\lambda)} E_{\gamma}E_{\gamma'-\gamma}\mathcal{T}_{\lambda}\mathcal{T}_{\lambda'-\lambda}f\\
&=\sum_{\substack{\lambda'\in\Lambda\\\gamma'\in\Gamma}}\langle g, E_{\gamma'-\gamma}\mathcal{T}_{\lambda'-\lambda}f\rangle  E_{\gamma}\mathcal{T}_{\lambda}E_{\gamma'-\gamma}\mathcal{T}_{\lambda'-\lambda}f\\
&=E_{\gamma}\mathcal{T}_{\lambda}T^{\ast}_{f:\Lambda,\Gamma}T_{f:\Lambda,\Gamma}(g),
\end{align*}
making use of the commutativity relation (\ref{E:commu relation}). This, in turn, implies that the operator $R_{\theta}^{-1}$ also commutes with $E_{\gamma}\mathcal{T}_{\lambda}$. Therefore,
\begin{align*}
T_{f:\Lambda,\Gamma}R_{\theta}^{-1}(g)
&=\{\langle g, E_{\gamma}\mathcal{T}_{\lambda}h_{\theta}\rangle\}_{(\lambda,\gamma)\in\Lambda\times\Gamma}
=T_{h_{\theta}:\Lambda,\Gamma}(g)
\end{align*}
and so
\begin{equation}\label{E:relation in Tstar f}
T^{\ast}_{f:\Lambda,\Gamma}T_{f:\Lambda,\Gamma}R_{\theta}^{-1}=T^{\ast}_{f:\Lambda,\Gamma}T_{h_{\theta}:\Lambda,\Gamma}.
\end{equation}
\par
Next, we define a linear functional $\psi$ on the space $\mathcal{B}(L^2(G))$ of bounded linear operators on $L^2(G)$ as follows. Let $E_i=(\alpha_i+\mathcal{S}_{\Gamma^{\perp}})\cap\mathcal{S}_{\Lambda}$ for $\alpha_i\in \Gamma^{\perp}$ and $i=1,2,\ldots,k$. Then, $E_i$'s are disjoint and by the hypothesis, we have $\bigcup\limits_{i=1}^k E_i=\mathcal{S}_{\Lambda}$. For $T\in\mathcal{B}(L^2(G))$, we define
\begin{equation*}
\psi(T):=\frac{1}{\mu_G(\mathcal{S}_{\Lambda})}\sum_{i=1}^k\langle T(\chi_{E_i}), \chi_{E_i}\rangle_{L^2(G)}.
\end{equation*}
Then,
\begin{align}
\psi(T^{\ast}_{f:\Lambda,\Gamma}T_{h_{\theta}:\Lambda,\Gamma})
&=\frac{1}{\mu_G(\mathcal{S}_{\Lambda})}\sum_{i=1}^k\langle T^{\ast}_{f:\Lambda,\Gamma}T_{h_{\theta}:\Lambda,\Gamma}(\chi_{E_i}), \chi_{E_i}\rangle_{L^2(G)}\nonumber\\
&=\frac{1}{\mu_G(\mathcal{S}_{\Lambda})}\sum_{i=1}^k\langle T_{h_{\theta}:\Lambda,\Gamma}(\chi_{E_i}), T_{f:\Lambda,\Gamma}(\chi_{E_i})\rangle_{l^2(\Lambda\times\Gamma)}\nonumber\\
&=\frac{1}{\mu_G(\mathcal{S}_{\Lambda})}\sum_{i=1}^k \sum_{\substack{\lambda'\in\Lambda\\\gamma'\in\Gamma}}(T_{h_{\theta}:\Lambda,\Gamma}(\chi_{E_i}))_{\lambda',\gamma'}\overline{(T_{f:\Lambda,\Gamma}(\chi_{E_i}))_{\lambda',\gamma'}}.\label{psi on T star Tf}
\end{align}
Now,
\begin{align*}
(T_{h_{\theta}:\Lambda,\Gamma}(\chi_{E_i}))_{\lambda',\gamma'}
&=\int_G \chi_{E_i}(x)\overline{\gamma'(x)h_{\theta}(x-\lambda')}d\mu_G(x)\\
&=\overline{\gamma'(\lambda')}\int_G \chi_{E_i}(x+\lambda')\overline{h_{\theta}(x)\gamma'(x)}d\mu_G(x)\\
&=\overline{\gamma'(\lambda')}\int_{E_i-\lambda'} \overline{h_{\theta}(x)\gamma'(x)}d\mu_G(x)\\
&=\int_{E_i-\alpha_i} \overline{h_{\theta}(x+\alpha_i-\lambda')\gamma'(x)}d\mu_G(x)\\
&=\int_{\mathcal{S}_{\Gamma^{\perp}}}\chi_{E_i}(x+\alpha_i) \overline{h_{\theta}(x+\alpha_i-\lambda')\gamma'(x)}d\mu_G(x)\\
&=\mu_G(\mathcal{S}_{\Gamma^{\perp}})\left(\chi_{E_i}(\cdot+\alpha_i) \overline{h_{\theta}(\cdot+\alpha_i-\lambda')}\right)^{\wedge}(\gamma').
\end{align*}
Using this in (\ref{psi on T star Tf}) and then applying the Parseval's identity, we obtain
\begin{align*}
&\psi(T^{\ast}_{f:\Lambda,\Gamma}T_{h_{\theta}:\Lambda,\Gamma})\\
&\,=\frac{\mu_G(\mathcal{S}_{\Gamma^{\perp}})^2}{\mu_G(\mathcal{S}_{\Lambda})}\sum_{i=1}^k \sum_{\substack{\lambda'\in\Lambda\\\gamma'\in\Gamma}}\left(\chi_{E_i}(\cdot+\alpha_i) \overline{h_{\theta}(\cdot+\alpha_i-\lambda')}\right)^{\wedge}(\gamma')\\
&\qquad\qquad\qquad\qquad\quad\times\overline{\left(\chi_{E_i}(\cdot+\alpha_i) \overline{f(\cdot+\alpha_i-\lambda')}\right)^{\wedge}(\gamma')}\\
&\,=\frac{\mu_G(\mathcal{S}_{\Gamma^{\perp}})}{\mu_G(\mathcal{S}_{\Lambda})}\sum_{i=1}^k \sum_{\lambda'\in\Lambda}\int\limits_{\mathcal{S}_{\Gamma^{\perp}}}\chi_{E_i}(x+\alpha_i) \overline{h_{\theta}(x+\alpha_i-\lambda')}f(x+\alpha_i-\lambda')d\mu_G(x)\\
&\,=\frac{\mu_G(\mathcal{S}_{\Gamma^{\perp}})}{\mu_G(\mathcal{S}_{\Lambda})} \sum_{\lambda'\in\Lambda}\sum_{i=1}^k\int_{E_i}f(x-\lambda')\overline{h_{\theta}(x-\lambda')}d\mu_G(x)\\
&\,=\frac{\mu_G(\mathcal{S}_{\Gamma^{\perp}})}{\mu_G(\mathcal{S}_{\Lambda})} \sum_{\lambda'\in\Lambda}\int_{\mathcal{S}_{\Lambda}-\lambda'}f(x)\overline{h_{\theta}(x)}d\mu_G(x)\\
&\,=\frac{\mu_G(\mathcal{S}_{\Gamma^{\perp}})}{\mu_G(\mathcal{S}_{\Lambda})}\langle f,h_{\theta}\rangle_{L^2(G)}.
\end{align*}
As $\mu_{\widehat{G}}(\mathcal{S}_{\Gamma})\mu_G(\mathcal{S}_{\Gamma^{\perp}})=1$, we further have
\begin{equation*}
\psi(T^{\ast}_{f:\Lambda,\Gamma}T_{h_{\theta}:\Lambda,\Gamma})=\frac{1}{\mu_G(\mathcal{S}_{\Lambda})\mu_{\widehat{G}}(\mathcal{S}_{\Gamma})}\langle f,h_{\theta}\rangle_{L^2(G)}.
\end{equation*}
It then follows from (\ref{E:relation in Tstar f}) that
\[\psi(T^{\ast}_{f:\Lambda,\Gamma}T_{f:\Lambda,\Gamma}R_{\theta}^{-1})=\frac{1}{\mu_G(\mathcal{S}_{\Lambda})\mu_{\widehat{G}}(\mathcal{S}_{\Gamma})}\langle f,h_{\theta}\rangle_{L^2(G)}.\]
We observe from the definition of $h_{\theta}$ that $\theta h_{\theta}+T^{\ast}_{f:\Lambda,\Gamma}T_{f:\Lambda,\Gamma}h_{\theta}=f$ and so
\begin{equation}\label{E:IP of f, h theta}
\langle f, h_{\theta}\rangle_{L^2(G)}=\theta\|h_{\theta}\|^2_{L^2(G)}+\|T_{f:\Lambda,\Gamma}h_{\theta}\|^2_{l^2(\Lambda\times\Gamma)}\geq \|T_{f:\Lambda,\Gamma}h_{\theta}\|^2_{l^2(\Lambda\times\Gamma)}\geq 0.
\end{equation}
Also,
\begin{align*}
\langle f, h_{\theta}\rangle_{L^2(G)}&=\langle h_{\theta},T^{\ast}_{f:\Lambda,\Gamma} e_{0_{\Lambda}0_{\Gamma}}\rangle,
\text{where } e_{\lambda'\gamma'}(\lambda,\gamma)=\delta_{\lambda\lambda'}\delta_{\gamma\gamma'}.
\end{align*}
Therefore,
\begin{align*}
\langle f, h_{\theta}\rangle_{L^2(G)}^2&=|\langle T_{f:\Lambda,\Gamma}h_{\theta}, e_{0_{\Lambda}0_{\Gamma}}\rangle|^2=|T_{f:\Lambda,\Gamma}h_{\theta}(0_{\Lambda},0_{\Gamma})|^2\leq \|T_{f:\Lambda,\Gamma}h_{\theta}\|^2_{l^2(\Lambda\times\Gamma)}\\
&\leq \langle f, h_{\theta}\rangle_{L^2(G)},
\end{align*}
by (\ref{E:IP of f, h theta}). This shows that $\langle f, h_{\theta}\rangle_{L^2(G)}\leq 1$, which further implies that
\[\psi(T^{\ast}_{f:\Lambda,\Gamma}T_{f:\Lambda,\Gamma}R_{\theta}^{-1})\leq \frac{1}{d(\Lambda\times\Gamma)}.\]
Taking the limit as $\theta$ tends to 0, we obtain
\begin{align*}
\frac{1}{\mu_G(\mathcal{S}_{\Lambda})}\sum_{i=1}^k\lim_{\theta\ra 0}\langle T^{\ast}_{f:\Lambda,\Gamma}T_{f:\Lambda,\Gamma}R_{\theta}^{-1}\chi_{E_i}, \chi_{E_i}\rangle_{L^2(G)}\leq \frac{1}{d(\Lambda\times\Gamma)}.
\end{align*}
It can be seen that the above limit actually exists. In fact, by Spectral theory, it is known that for any $g,\,h\in L^2(G)$,
$$\langle T^{\ast}_{f:\Lambda,\Gamma}T_{f:\Lambda,\Gamma}R_{\theta}^{-1}g, h\rangle_{L^2(G)}\ra\langle Pg,h\rangle_{L^2(G)}\text{ as }\theta\ra 0,$$ where $P$ denotes the projection operator on the closure of the range of $T^{\ast}_{f:\Lambda,\Gamma}$. By hypothesis, the Gabor system $\{E_{\gamma}\mathcal{T}_{\lambda}f\}_{(\lambda,\gamma)\in\Lambda\times\Gamma}$ is complete and so the range of $T^{\ast}_{f:\Lambda,\Gamma}$ is dense in $L^2(G)$. In other words, the operator $P$ turns out to be the identity operator in this case. Therefore,
\[\frac{1}{\mu_G(\mathcal{S}_{\Lambda})}\sum_{i=1}^k\mu_G(E_i)\leq \frac{1}{d(\Lambda\times\Gamma)}.\]
As $\mathcal{S}_{\Lambda}$ is the disjoint union of these $E_i$'s, we get $d(\Lambda\times\Gamma)\leq 1$, thereby proving the theorem.
\end{proof}

\begin{remark}
The above result can alternatively be obtained as an application of Theorem 7.4 in \cite{romero2022density}, which involves the theory of representations. In this direction, we need to first consider A. Weil's abstract Heisenberg group $\mathbb{H}(G)=G\times\widehat{G}\times\mathbb{T}$, wherein the group operation is given by 
\[(x,\nu,\theta)\cdot(x',\nu',\theta')=(x+x',\nu+\nu',\theta\theta'\nu'(x)).\] By the Stone-von Neumann theorem, it is known that every infinite dimensional irreducible unitary representation of $\mathbb{H}(G)$ is unitarily equivalent to the representations $\rho_j$, with $j$ a non-zero integer, on the representation space $L^2(G)$ and given by 
\[\rho_j(x,\nu,\theta)f(u)=\theta^j(\nu(u))^jf(x+u),\quad f\in L^2(G).\]
We refer to Folland \cite{Foll} for further details. Being closely connected with the above action, the map $\pi$ from $G\times\widehat{G}$ into the group of unitary operators acting on $L^2(G)$, defined by
\[\pi(x,\nu)f(u)=\nu(u)f(u-x),\quad f\in L^2(G),\]
turns out to be an irreducible $\sigma-$representation with cocyle $\sigma((x,\nu),(x',\nu'))$ $:=\nu'(-x).$ It can be seen that the Fourier-Wigner transform $V$, given by $V(f,g)(x,\nu)=\langle \pi(x,\nu)f,g\rangle$, for $f,g\in L^2(G)$ and $(x,\nu)\in G\times\widehat{G}$, satisfies
$V(f,g)\in L^2(G\times\widehat{G})$ and 
\[\int_{G}\int_{\widehat{G}}V(f_1,g_1)(x,\nu)\overline{V(f_2,g_2)(x,\nu)}d\mu_{\widehat{G}}(\nu)d\mu_{G}(x)=\langle f_1, f_2\rangle\overline{\langle g_1, g_2\rangle},\] as in the classical case. So, $\pi$ is also square-integrable and has 1 as its formal dimension. Now, taking the lattice $\Lambda\times\Gamma$ in the group $G\times\widehat{G}$, the orbit of the restriction of $(\pi,L^2(G))$ to $\Lambda\times\Gamma$ is the Gabor system $\{E_{\gamma}\mathcal{T}_{\lambda}f\}_{(\lambda,\gamma)\in\Lambda\times\Gamma}$, for some $f\in L^2(G)$. The result then follows by applying Theorem 7.4 in \cite{romero2022density} to the lattice $\Lambda\times\Gamma$ and the $\sigma-$representation $\pi$ of $G\times\widehat{G}$.
\end{remark}

We provide both these proofs as both of them are elegant in their own respect. We would also like to remark that though Bekka in \cite{bekka2004square} has proved this result for $\mathbb{R}^n$, using the classical Heisenberg group, that proof cannot be directly lifted to the abstract setting, which doesn't have a nilpotent Lie group structure.  

\section{Adjoint Gabor system as R-dual}\label{S:Adjoint Gabor system as R-dual}
The following theorem shows that the adjoint Gabor system $\{E_{\beta}\mathcal{T}_\alpha f\}_{(\alpha,\beta)\in\Gamma^{\perp}\times\Lambda^{\perp}}$ is an R-dual, when the given the Gabor frame $\{E_{\gamma}\mathcal{T}_{\lambda}f\}_{(\lambda,\gamma)\in\Lambda\times\Gamma}$ is tight and $d(\Lambda\times\Gamma)<1$. 
\begin{theorem}\label{T:adj for 1-tight}
Let the lattices $\Lambda$ of $G$ and $\Gamma$ of $\widehat{G}$ be such that both $\Gamma$ and $\Gamma^{\perp}$ contain an element of infinite order, there exist a finite number, say $k$, of elements $\alpha_i\in \Lambda$ such that $\bigcup\limits_{i=1}^k (\alpha_i+\mathcal{S}_{\Lambda})\cap\mathcal{S}_{\Gamma^{\perp}}=\mathcal{S}_{\Gamma^{\perp}}$ and $d(\Lambda\times\Gamma)<1$. For $f\in L^2(G)$, suppose $\{E_{\gamma}\mathcal{T}_{\lambda}f\}_{(\lambda,\gamma)\in\Lambda\times\Gamma}$ is a tight frame for $L^2(G)$. Then, the adjoint Gabor system $\{E_{\beta}\mathcal{T}_\alpha f\}_{(\alpha,\beta)\in\Gamma^{\perp}\times\Lambda^{\perp}}$ is an R-dual of the Gabor frame $\{E_{\gamma}\mathcal{T}_{\lambda}f\}_{(\lambda,\gamma)\in\Lambda\times\Gamma}$.
\end{theorem}

\begin{proof}
We may assume that the given Gabor system is a Parseval frame without any loss of generality. For $\lambda'\in\Lambda$ and $\gamma'\in\Gamma$, consider the operators $U_{\lambda'}$ and $V_{\gamma'}$ on $l^2(\Lambda\times\Gamma)$ defined as follows:
For $\{c_{\lambda,\gamma}\}\in l^2(\Lambda\times\Gamma)$,
\begin{equation*}
U_{\lambda'}(\{c_{\lambda,\gamma}\}):=\{\overline{\gamma(\lambda')}c_{\lambda-\lambda',\gamma}\}\quad\text{ and }\quad
V_{\gamma'}(\{c_{\lambda,\gamma}\}):=\{c_{\lambda,\gamma-\gamma'}\}.
\end{equation*}
Clearly, both these operators are isometries on $l^2(\Lambda\times\Gamma)$. These operators are needed to establish  intertwining relations of the analysis operator $T_{f:\Lambda,\Gamma}$ with the translation and modulation operators. In fact, for $\lambda'\in\Lambda$ and $g\in L^2(G)$,
\begin{align*}
T_{f:\Lambda,\Gamma}\mathcal{T}_{\lambda'}g
&=\{\langle \mathcal{T}_{\lambda'}g, E_{\gamma}\mathcal{T}_{\lambda}f\rangle\}_{(\lambda,\gamma)\in\Lambda\times\Gamma}
=\{\overline{\gamma(\lambda')}\langle g, E_{\gamma}\mathcal{T}_{\lambda-\lambda'}f\rangle\}_{(\lambda,\gamma)\in\Lambda\times\Gamma}\\
&=U_{\lambda'}(\{\langle g, E_{\gamma}\mathcal{T}_{\lambda}f\rangle\}_{(\lambda,\gamma)\in\Lambda\times\Gamma})
=U_{\lambda'}T_{f:\Lambda,\Gamma}g,
\end{align*}
using (\ref{E:commu relation}). Also, for $\gamma'\in\Gamma$,
\begin{align*}
T_{f:\Lambda,\Gamma}E_{\gamma'}g
&=\{\langle E_{\gamma'}g, E_{\gamma}\mathcal{T}_{\lambda}f\rangle\}_{(\lambda,\gamma)\in\Lambda\times\Gamma}
=\{\langle g, E_{\gamma-\gamma'}\mathcal{T}_{\lambda}f\rangle\}_{(\lambda,\gamma)\in\Lambda\times\Gamma}\\
&=V_{\gamma'}(\{\langle g, E_{\gamma}\mathcal{T}_{\lambda}f\rangle\}_{(\lambda,\gamma)\in\Lambda\times\Gamma})
=V_{\gamma'}T_{f:\Lambda,\Gamma}g.
\end{align*}
This shows that both $(\mathcal{T}_{\lambda'},U_{\lambda'})$ and $(E_{\gamma'},V_{\gamma'})$ are intertwining pairs of operators for $\{E_{\gamma}\mathcal{T}_{\lambda}f\}_{(\lambda,\gamma)\in\Lambda\times\Gamma}$.
\par
Now, for $\gamma'\in \Gamma$ of infinite order and $a\in\mathbb{C}$, let $c\in l^2(\Lambda\times\Gamma)$ be such that $V_{\gamma'}c=ac$. So, $c_{\lambda,\gamma-\gamma'}=ac_{\lambda,\gamma},\,\forall\,(\lambda,\gamma)\in\Lambda\times\Gamma$. Then, by taking $[\gamma']$ to denote the cyclic subgroup generated by $\gamma'$, we have
\begin{align*}
\|c\|_{l^2(\Lambda\times\Gamma)}^2
&=\sum_{\lambda\in\Lambda}\int_{\Gamma/[\gamma']}\sum_{k\in\mathbb{Z}}|c_{\lambda,\gamma+k\gamma'}|^2
d\mu_{\Gamma/[\gamma']}(\gamma+[\gamma'])\\
&=\sum_{\lambda\in\Lambda}\sum_{\gamma\in\mathcal{S}_{[\gamma']}}\sum_{k\in\mathbb{Z}}|c_{\lambda,\gamma-k\gamma'}|^2\\
&=\sum_{\lambda\in\Lambda}\sum_{\gamma\in\mathcal{S}_{[\gamma']}}|c_{\lambda,\gamma}|^2\sum_{k\in\mathbb{Z}}|a|^{2k}.
\end{align*}
This forces $a$ to be $0$, as $\|c\|_{l^2(\Lambda\times\Gamma)}^2<\infty$ and so $c$ also turns out to be the zero sequence. In other words, $\text{ker}(V_{\gamma'}-aI)=\{0\}$. Further, as $V_{\gamma'}^{\ast}=V_{-\gamma'}$, by (\ref{E:excess and synth ker}) and Theorem \ref{T:intertwining}, dim(ker $T_{f:\Lambda,\Gamma}^{\ast}$), which is equal to the excess of the frame $\{E_{\gamma}\mathcal{T}_{\lambda}f\}_{(\lambda,\gamma)\in\Lambda\times\Gamma}$, is either zero or infinite. By the condition $d(\Lambda\times\Gamma)<1$ in the hypothesis, we have that $\{E_{\gamma}\mathcal{T}_{\lambda}f\}_{(\lambda,\gamma)\in\Lambda\times\Gamma}$ is not a Riesz basis, by Theorem 5.7 in \cite{jakobsen2016density}. This implies that dim(ker$(T_{f:\Lambda,\Gamma}^{\ast})$) is non-zero and hence infinity.
\par
On the other hand, $(\mathcal{T}_{\alpha'},U_{\alpha'})$ is an intertwining pair of operators for the Bessel sequence $\{E_{\beta}\mathcal{T}_\alpha f\}_{(\alpha,\beta)\in\Gamma^{\perp}\times\Lambda^{\perp}}$ in $L^2(G)$, where for $\alpha'\in\Gamma^{\perp}$, $U_{\alpha'}$ is an operator on $l^2(\Gamma^{\perp}\times\Lambda^{\perp})$ defined analogously. Now, let $\alpha'\in\Gamma^{\perp}$ be of infinite order and $g\in \text{ker}(\mathcal{T}_{\alpha'}-aI)$, $a\in\mathbb{C}$. In other words, 
\begin{equation}\label{E:pt spec eqn}
g(x-\alpha')=ag(x)\text{ a.e. }x\in G.
\end{equation}
If $a=0$, then by (\ref{E:pt spec eqn}), $g=0$ in $L^2(G)$ and so $\text{ker}(\mathcal{T}_{\alpha'}-aI)=\{0\}$. Suppose $a\neq 0$. Then, for $\tilde{\alpha}\in\Gamma^{\perp}$,
\begin{align}
\int\limits_{[\alpha]+\mathcal{S}_{\Gamma^{\perp}}+\tilde{\alpha}}|g(x)|^2d\mu_G(x)
&=\sum_{k\in\mathbb{Z}}\int\limits_{\mathcal{S}_{\Gamma^{\perp}}+\tilde{\alpha}}|g(k\alpha+x)|^2d\mu_{G}(x)\nonumber\\
&=\sum_{k\in\mathbb{Z}}\int\limits_{\mathcal{S}_{\Gamma^{\perp}}+\tilde{\alpha}}|a|^{-2k}|g(x)|^2d\mu_{G}(x)\nonumber\\
&=\sum_{k\in\mathbb{Z}}|a|^{-2k}\int\limits_{\mathcal{S}_{\Gamma^{\perp}}+\tilde{\alpha}}|g(x)|^2d\mu_{G}(x),\label{E:pt spec eqn 2}
\end{align}
using (\ref{E:pt spec eqn}). As $g\in L^2(G)$ and $a\neq 0$, from (\ref{E:pt spec eqn 2}), we get $g=0$ a.e. on $\mathcal{S}_{\Gamma^{\perp}}+\tilde{\alpha}$. This in turn implies that $g=0$ in $L^2(G)$ and hence, $\text{ker}(\mathcal{T}_{\alpha'}-aI)=\{0\}$. By Theorem \ref{T:intertwining}, the deficit of $\{E_{\beta}\mathcal{T}_\alpha f\}_{(\alpha,\beta)\in\Gamma^{\perp}\times\Lambda^{\perp}}$ is either zero or infinity. We have $d(\Gamma^{\perp}\times\Lambda^{\perp})>1$ and so by Theorem \ref{T:GF completeness}, the Bessel sequence $\{E_{\beta}\mathcal{T}_\alpha f\}_{(\alpha,\beta)\in\Gamma^{\perp}\times\Lambda^{\perp}}$ is not complete in $L^2(G)$, thereby having non-zero deficit. Thus, the infinite deficit of this sequence results in dim((span$\{E_{\beta}\mathcal{T}_\alpha f:(\alpha,\beta)\in\Gamma^{\perp}\times\Lambda^{\perp}\})^{\perp})=\infty$.
\par
Now, let $\{e_{\alpha\beta}\}_{(\alpha,\beta)\in\Gamma^{\perp}\times\Lambda^{\perp}}$ and $\{\tilde{e}_{\lambda\gamma}\}_{(\lambda,\gamma)\in\Lambda\times\Gamma}$ be orthonormal bases for $L^2(G)$. We shall define for $\alpha\in\Gamma^{\perp},\,\beta\in\Lambda^{\perp}$,
\[w_{\alpha\beta}:=\sum_{\substack{\lambda\in\Lambda\\\gamma\in\Gamma}}\langle E_{\gamma}\mathcal{T}_{\lambda}f, e_{\alpha\beta}\rangle_{L^2(G)}\tilde{e}_{\lambda\gamma}.\]
Then, for $\alpha,\,\alpha'\in\Gamma^{\perp},\,\beta,\,\beta'\in\Lambda^{\perp}$,
\begin{align*}
\langle w_{\alpha\beta},w_{\alpha'\beta'}\rangle
&=\sum_{\lambda,\gamma}\sum_{\lambda',\gamma'}\langle E_{\gamma}\mathcal{T}_{\lambda}f, e_{\alpha\beta}\rangle_{L^2(G)} \langle e_{\alpha'\beta'}, E_{\gamma'}\mathcal{T}_{\lambda'}f\rangle_{L^2(G)} \langle \tilde{e}_{\lambda\gamma}, \tilde{e}_{\lambda'\gamma'}\rangle_{L^2(G)}\\
&=\sum_{\lambda,\gamma}\langle E_{\gamma}\mathcal{T}_{\lambda}f, e_{\alpha\beta}\rangle_{L^2(G)} \langle e_{\alpha'\beta'}, E_{\gamma}\mathcal{T}_{\lambda}f\rangle_{L^2(G)}\\
&=\left\langle \sum_{\lambda,\gamma}\langle e_{\alpha'\beta'}, E_{\gamma}\mathcal{T}_{\lambda}f\rangle_{L^2(G)}E_{\gamma}\mathcal{T}_{\lambda}f, e_{\alpha\beta}\right\rangle_{L^2(G)}\\
&=\langle S_{f:\Lambda,\Gamma}e_{\alpha'\beta'}, e_{\alpha\beta}\rangle_{L^2(G)}\\
&=\delta_{\alpha\alpha'}\delta_{\beta\beta'},
\end{align*}
as $\{E_{\gamma}\mathcal{T}_{\lambda}f\}_{(\lambda,\gamma)\in\Lambda\times\Gamma}$ is a Parseval frame. Therefore, $\{w_{\alpha\beta}\}$ is an orthonormal basis for its closed linear span. We shall now prove that dim(($\text{span }\{w_{\alpha\beta}:(\alpha,\beta)\in\Gamma^{\perp}\times\Lambda^{\perp}\})^{\perp})$ is also infinite. For $g\in L^2(G)$, we have
\begin{align*}
\langle g, w_{\alpha\beta}\rangle&=\left\langle g,\sum_{\substack{\lambda\in\Lambda\\\gamma\in\Gamma}}\langle E_{\gamma}\mathcal{T}_{\lambda}f, e_{\alpha\beta}\rangle_{L^2(G)}\tilde{e}_{\lambda\gamma}\right\rangle\\
&=\sum_{\substack{\lambda\in\Lambda\\\gamma\in\Gamma}}\langle e_{\alpha\beta}, E_{\gamma}\mathcal{T}_{\lambda}f\rangle_{L^2(G)}\langle g,\tilde{e}_{\lambda\gamma}\rangle\\
&=\left\langle e_{\alpha\beta},\sum_{\substack{\lambda\in\Lambda\\\gamma\in\Gamma}}\langle \tilde{e}_{\lambda\gamma}, g\rangle E_{\gamma}\mathcal{T}_{\lambda}f\right\rangle.
\end{align*}
So, we may define an operator $R: (\text{span }\{w_{\alpha\beta}:(\alpha,\beta)\in\Gamma^{\perp}\times\Lambda^{\perp}\})^{\perp}\ra \text{ker }T_{f:\Lambda,\Gamma}^{\ast}$ by $Rg:=\{\langle \tilde{e}_{\lambda\gamma}, g\rangle\}_{(\lambda,\gamma)\in\Lambda\times\Gamma}$, which is a conjugate-linear surjective isometry. We had shown earlier that the dimension of the kernel of $T_{f:\Lambda,\Gamma}^{\ast}$ is infinity. Therefore, dim(($\text{span }\{w_{\alpha\beta}:(\alpha,\beta)\in\Gamma^{\perp}\times\Lambda^{\perp}\})^{\perp})=\infty$. We may take $\{\varphi_{\mu\nu}\}$ and $\{\psi_{\mu\nu}\}$ to be orthonormal bases of $(\text{span }\{w_{\alpha\beta}:(\alpha,\beta)\in\Gamma^{\perp}\times\Lambda^{\perp}\})^{\perp}$ and (span $\{E_{\beta}\mathcal{T}_\alpha f:(\alpha,\beta)\in\Gamma^{\perp}\times\Lambda^{\perp}\})^{\perp}$ respectively. Moreover, by Theorem \ref{T:dual orthog}, $\left\{E_{\beta}\mathcal{T}_\alpha f\right\}$ is an orthonormal basis for its closed linear span. 
\par
Now, as $\{w_{\alpha\beta}\}$ is an orthonormal basis for $\text{span }\{w_{\alpha\beta}:(\alpha,\beta)\in\Gamma^{\perp}\times\Lambda^{\perp}\}$ and $\{\varphi_{\mu\nu}\}$ is one for its orthogonal complement, we define a linear operator $\mathcal{U}: L^2(G)\ra L^2(G)$ by
\[\mathcal{U}(w_{\alpha\beta})=E_{\beta}\mathcal{T}_\alpha f\text{ and }\mathcal{U}(\varphi_{\mu\nu})=\psi_{\mu\nu}.\]  This operator $\mathcal{U}$ takes an orthonormal basis of $L^2(G)$ to another orthornormal basis and is therefore a unitary operator on $L^2(G)$ such that
\[E_{\beta}\mathcal{T}_\alpha f=\mathcal{U}(w_{\alpha\beta})=\sum_{\substack{\lambda\in\Lambda\\\gamma\in\Gamma}}\langle E_{\gamma}\mathcal{T}_{\lambda}f, e_{\alpha\beta}\rangle_{L^2(G)}\mathcal{U}(\tilde{e}_{\lambda\gamma}).\]
Thus, the adjoint Gabor system, $\{E_{\beta}\mathcal{T}_\alpha f\}_{(\alpha,\beta)\in\Gamma^{\perp}\times\Lambda^{\perp}}$ is an R-dual of $\{E_{\gamma}\mathcal{T}_{\lambda}f\}_{(\lambda,\gamma)\in\Lambda\times\Gamma}$ with respect to the orthonormal bases $\{e_{\alpha\beta}\}_{(\alpha,\beta)\in\Gamma^{\perp}\times\Lambda^{\perp}}$ and $\{\mathcal{U}(\tilde{e}_{\lambda\gamma})\}_{(\lambda,\gamma)\in\Lambda\times\Gamma}$.
\end{proof}
We now consider the special case when $\Gamma=\Lambda^{\perp}$. The proposition below gives an explicit construction of an orthonormal basis with Gabor structure for $L^2(G)$, which will be needed for the R-duality (Theorem \ref{T:R duality for Lambda and Lambda perp}) of the Gabor system associated with the lattice $\Lambda\times\Lambda^{\perp}$. Though the proof of the construction is straightforward and follows using standard arguments, we provide the same for the sake of completeness.
\begin{prop}\label{P:Exis of OB}
The collection $
\left\{\frac{1}{\sqrt{\mu_G(\mathcal{S}_{\Lambda})}}E_{\gamma}\mathcal{T}_{\lambda}\chi_{\mathcal{S}_{\Lambda}}\right\}_{(\lambda,\gamma)\in\Lambda\times\Lambda^{\perp}}$ is an orthonormal basis for $L^2(G)$.
\end{prop}

\begin{proof}
For $(\lambda,\gamma)\in\Lambda\times\Lambda^{\perp},$ clearly $\left\|\frac{1}{\sqrt{\mu_G(\mathcal{S}_{\Lambda})}}E_{\gamma}\mathcal{T}_{\lambda}\chi_{\mathcal{S}_{\Lambda}}\right\|_{L^2(G)}=1$ and for $(\lambda,\gamma),\,(\lambda',\gamma')\in\Lambda\times\Lambda^{\perp}$ with $\lambda\neq\lambda'$, we have
\begin{align*}
&\left\langle \frac{1}{\sqrt{\mu_G(\mathcal{S}_{\Lambda})}}E_{\gamma}\mathcal{T}_{\lambda}\chi_{\mathcal{S}_{\Lambda}},\frac{1}{\sqrt{\mu_G(\mathcal{S}_{\Lambda})}}E_{\gamma'}\mathcal{T}_{\lambda'}\chi_{\mathcal{S}_{\Lambda}}\right
\rangle_{L^2(G)}\\
&\;\quad\qquad=\frac{1}{\mu_G(\mathcal{S}_{\Lambda})}\int_G\gamma(x)\overline{\gamma'(x)}\chi_{\mathcal{S}_{\Lambda}}(x-\lambda)\overline{\chi_{\mathcal{S}_{\Lambda}}(x-\lambda')}d\mu_G(x).
\end{align*}
As $\lambda\neq\lambda'$, both $x-\lambda$ and $x-\lambda'$, for any $x\in G$, cannot belong to $\mathcal{S}_{\Lambda}$ at the same time and so, the above inner product is zero. Now, let $\lambda=\lambda'$ but $\gamma\neq\gamma'$. Then, there exists $x_0\in \mathcal{S}_{\Lambda}$ such that $\gamma(x_0)\neq\gamma'(x_0)$. Consider,
\begin{align*}
&\left\langle \frac{1}{\sqrt{\mu_G(\mathcal{S}_{\Lambda})}}E_{\gamma}\mathcal{T}_{\lambda}\chi_{\mathcal{S}_{\Lambda}},\frac{1}{\sqrt{\mu_G(\mathcal{S}_{\Lambda})}}E_{\gamma'}\mathcal{T}_{\lambda'}\chi_{\mathcal{S}_{\Lambda}}\right
\rangle_{L^2(G)}\\
&\;\quad\qquad=\frac{1}{\mu_G(\mathcal{S}_{\Lambda})}\int_G\gamma(x)\overline{\gamma'(x)}\chi_{\mathcal{S}_{\Lambda}}(x-\lambda)d\mu_G(x)\\
&\;\quad\qquad=\frac{1}{\mu_G(\mathcal{S}_{\Lambda})}\int_G(\gamma-\gamma')(u)\chi_{\mathcal{S}_{\Lambda}}(u)d\mu_G(u).
\end{align*}
Now,
\begin{align*}
&\int_G(\gamma-\gamma')(u)\chi_{\mathcal{S}_{\Lambda}}(u)d\mu_G(u)
=(\gamma-\gamma')(x_0)\int_G (\gamma-\gamma')(u)\chi_{\mathcal{S}_{\Lambda}}(u)d\mu_G(u).
\end{align*}
This forces $\int_G(\gamma-\gamma')(u)\chi_{\mathcal{S}_{\Lambda}}(u)d\mu_G(u)$ to be zero and hence the inner product is zero again in this case. Thus, we have shown that $
\left\{\frac{1}{\sqrt{\mu_G(\mathcal{S}_{\Lambda})}}E_{\gamma}\mathcal{T}_{\lambda}\chi_{\mathcal{S}_{\Lambda}}\right\}_{(\lambda,\gamma)\in\Lambda\times\Lambda^{\perp}}$ is orthonormal in $L^2(G)$.
\par
In order to show completeness, consider $f\in L^2(G)$ such that $$\left\langle f,\frac{1}{\sqrt{\mu_G(\mathcal{S}_{\Lambda})}}E_{\gamma}\mathcal{T}_{\lambda}\chi_{\mathcal{S}_{\Lambda}}\right
\rangle_{L^2(G)}=0,\quad\forall\,(\lambda,\gamma)\in\Lambda\times\Lambda^{\perp}.$$
Let $\lambda\in \Lambda$. Then, for any $\gamma\in\Lambda^{\perp}$,
\begin{align*}
0&=\frac{1}{\sqrt{\mu_G(\mathcal{S}_{\Lambda})}}\int_G f(x)\overline{\gamma(x)}\chi_{\mathcal{S}_{\Lambda}}(x-\lambda)d\mu_G(x)\\
&=\frac{1}{\sqrt{\mu_G(\mathcal{S}_{\Lambda})}}\int_{\mathcal{S}_{\Lambda}} \mathcal{T}_{-\lambda}f(x)\overline{\gamma(x)}d\mu_G(x).
\end{align*}
In other words, $(\mathcal{T}_{-\lambda}f)^{\wedge}(\gamma)=0$, $\forall\,\gamma\in\Lambda^{\perp}$. By Plancherel theorem, $\mathcal{T}_{-\lambda}f=0$ in $L^2(\mathcal{S}_{\Lambda})$. So, $f(x)=0$ a.e. $x\in\mathcal{S}_{\Lambda}+\lambda$. As $G=\bigcup\limits_{\lambda\in\Lambda}(\mathcal{S}_{\Lambda}+\lambda)$ and $\Lambda$ is countable, we get $f=0$ almost everywhere on $G$, thereby proving the proposition.
\end{proof}

\begin{theorem}\label{T:R duality for Lambda and Lambda perp}
There exist orthonormal bases in $L^2(G)$ such that the adjoint Gabor system $\{E_{\beta}\mathcal{T}_\alpha f\}_{(\alpha,\beta)\in\Lambda\times\Lambda^{\perp}}$ is an R-dual of the Gabor system $\{E_{\gamma}\mathcal{T}_{\lambda}f\}_{(\lambda,\gamma)\in\Lambda\times\Lambda^{\perp}}$, for every $f\in L^2(G)$.  
\end{theorem}

\begin{proof}
By Proposition \ref{P:Exis of OB}, $\left\{\frac{1}{\sqrt{\mu_G(\mathcal{S}_{\Lambda})}}E_{\gamma}\mathcal{T}_{\lambda}\chi_{\mathcal{S}_{\Lambda}}\right\}_{(\lambda,\gamma)\in\Lambda\times\Lambda^{\perp}}$ is an orthonormal basis for $L^2(G)$. For $(\lambda,\gamma),\,(\alpha,\beta)\in\Lambda\times\Lambda^{\perp}$, we have
\begin{align*}
\langle E_{\gamma}\mathcal{T}_{\lambda}f,E_{\beta}\mathcal{T}_{\alpha}\chi_{\mathcal{S}_{\Lambda}}\rangle_{L^2(G)}
&=\langle E_{-\beta}\mathcal{T}_{-\alpha}f,E_{-\gamma}\mathcal{T}_{-\lambda}\chi_{\mathcal{S}_{\Lambda}}\rangle_{L^2(G)}.
\end{align*}
Consequently,
\begin{align*}
&\sum_{\substack{\lambda\in\Lambda\\\gamma\in\Lambda^{\perp}}}\left|\left\langle E_{\gamma}\mathcal{T}_{\lambda}f, \frac{1}{\sqrt{\mu_G(\mathcal{S}_{\Lambda})}}E_{-\beta}\mathcal{T}_{-\alpha}\chi_{\mathcal{S}_{\Lambda}}\right
\rangle_{L^2(G)}\right|^2\\
&\qquad\qquad=\sum_{\substack{\lambda\in\Lambda\\\gamma\in\Lambda^{\perp}}}\left|\left\langle E_{\beta}\mathcal{T}_{\alpha}f, \frac{1}{\sqrt{\mu_G(\mathcal{S}_{\Lambda})}}E_{-\gamma}\mathcal{T}_{-\lambda}\chi_{\mathcal{S}_{\Lambda}}\right
\rangle_{L^2(G)}\right|^2<\infty,
\end{align*}
by Proposition \ref{P:Exis of OB}. Therefore, the R-dual of $\{E_{\gamma}\mathcal{T}_{\lambda}f\}_{(\lambda,\gamma)\in\Lambda\times\Lambda^{\perp}}$ with respect to the orthonormal bases $\left\{\frac{1}{\sqrt{\mu_G(\mathcal{S}_{\Lambda})}}E_{\beta}\mathcal{T}_{\alpha}\chi_{\mathcal{S}_{\Lambda}}\right\}$ and $\left\{\frac{1}{\sqrt{\mu_G(\mathcal{S}_{\Lambda})}}E_{\gamma}\mathcal{T}_{\lambda}\chi_{\mathcal{S}_{\Lambda}}\right\}$ is given by 
\begin{align*}
&\sum_{\substack{\lambda\in\Lambda\\\gamma\in\Lambda^{\perp}}}\left\langle E_{\gamma}\mathcal{T}_{\lambda}f, \frac{1}{\sqrt{\mu_G(\mathcal{S}_{\Lambda})}}E_{-\beta}\mathcal{T}_{-\alpha}\chi_{\mathcal{S}_{\Lambda}}\right
\rangle_{L^2(G)}\frac{1}{\sqrt{\mu_G(\mathcal{S}_{\Lambda})}}E_{-\gamma}\mathcal{T}_{-\lambda}\chi_{\mathcal{S}_{\Lambda}}\\
&\;=\sum_{\substack{\lambda\in\Lambda\\\gamma\in\Lambda^{\perp}}}\left\langle E_{\beta}\mathcal{T}_{\alpha}f, \frac{1}{\sqrt{\mu_G(\mathcal{S}_{\Lambda})}}E_{-\gamma}\mathcal{T}_{-\lambda}\chi_{\mathcal{S}_{\Lambda}}\right
\rangle_{L^2(G)}\frac{1}{\sqrt{\mu_G(\mathcal{S}_{\Lambda})}}E_{-\gamma}\mathcal{T}_{-\lambda}\chi_{\mathcal{S}_{\Lambda}}\\
&\;= E_{\beta}\mathcal{T}_{\alpha}f,
\end{align*}
thereby proving the theorem.
\end{proof}

\begin{remark}
One can easily see that the above theorem holds true even for the finite abelian group setting, by following the lines of proofs of Proposition \ref{P:Exis of OB} and Theorem \ref{T:R duality for Lambda and Lambda perp} with the integrals being replaced by summation.  
\end{remark}

\noindent {\bf Acknowledgements:}
The first named author, S. Arati, would like to thank the National Board for Higher Mathematics, Department of Atomic Energy(Government of India) for the funding.

\end{document}